\documentclass[reqno, 10pt]{amsart}

\usepackage{amsmath,amssymb,amsthm}

\usepackage{tikz}
\usepackage{caption}
\usepackage{graphicx}
\usepackage{graphics}
\usepackage{tikz}	
\usepackage{MnSymbol}
\usetikzlibrary{cd,positioning, shapes}	
\usepackage{comment}
\usepackage{quiver}
\usepackage{mdframed}
\usepackage{empheq}
\usepackage[hidelinks]{hyperref}
\usepackage{longtable}
\usepackage{float}
\usepackage{standalone}
\usepackage[normalem]{ulem}

\newtheorem*{thm}{Theorem}
\newtheorem{theorem}{Theorem}[section]
\newtheorem{proposition}[theorem]{Proposition}
\newtheorem{question}[theorem]{Question}

\newtheorem{conjecture}[theorem]{Conjecture}

\theoremstyle{definition}

\newtheorem{remark}{Remark}[section]

\newcommand{\RP}{\mathbb{R}P}
\newcommand{\bbS}{\mathbb{S}}
\DeclareMathOperator{\str}{st}

\DeclareMathOperator{\conv}{conv}
\newcommand{\e}{\mathbf{e}}

\begin{document}

\title[]{An Efficient Triangulation of $\mathbb{R}P^5$}

\author[D. Guyer]{Dan Guyer}
\address[D. Guyer]{Department of Mathematics, University of Washington, Seattle, WA.}
\email{dguyer@uw.edu}

\author[S. Steinerberger]{Stefan Steinerberger}
\address[S. Steinerberger]{Department of Mathematics and Department of Applied Mathematics, University of Washington, Seattle, WA.}
\email{steinerb@uw.edu}

\author[Y. Yang]{Yirong Yang}
\address[Y. Yang]{School of Mathematics, Georgia Institute of Technology, Atlanta, GA.}
\email{yyang917@gatech.edu}

\begin{abstract}
    We present a $6$-dimensional centrally symmetric simplicial polytope for which the antipodal quotient of its boundary forms a $24$-vertex triangulation of the $5$-dimensional real projective space. This $6$-polytope is highly symmetric with an automorphism group of order $192$, and is of independent interest. We conjecture that our construction uses the fewest number of vertices among all triangulations of $\mathbb{R}P^5$. Our method also produces two triangulations of $\mathbb{R}P^6$ on $45$ and $49$ vertices; both improve the previously best known construction in dimension $6$ that used 53 vertices.
\end{abstract}

\maketitle

\section{Introduction and Results}
\subsection{Introduction} Simplicial complexes are a foundational way to combinatorially represent a topological space. It is a classical and hard problem in piecewise-linear topology to find \emph{vertex-minimal triangulations} of a given manifold. When the manifold is the sphere $\bbS^d$, one can take the boundary of the $(d+1)$-dimensional simplex, which has $d+2$ vertices. Additionally, the answer is known for PL triangulations of $\bbS^{d-1}\times \bbS^1$ and the twisted bundle $\bbS^{d-1} \dtimes \bbS^1$~\cite{Kuhnel1986,Bagchi2008,Chestnut_2008}. However, when the homology or homotopy groups of the manifold become more complicated, the number of vertices needed in its triangulation tends to increase. For a survey on efficient triangulations, we recommend Adiprasito--Benedetti~\cite{AdiprasitoBenedetti2025}.
\\

We study this problem for another central
topological space, the $d$-dimensional real projective space $\RP^d$. Even here, the bounds for the appropriate number of vertices remain very far apart. The best known lower bound comes from Arnoux and Marin \cite{Arnoux_1991} using the cup length of the real projective space (see Remark~\ref{rem:frick}).

\begin{proposition}[Arnoux--Marin, 1991]
    If $d\geq 3$, then any triangulation of real projective space $\RP^d$ must have at least $(d+2)(d+1)/2 + 1$ vertices. 
\end{proposition}
\noindent
Other methods to provide lower bounds include using topological estimates from the LS-category~\cite{Govc2020} and commutative algebra techniques~\cite{Novik2009,Murai2015}. \\

Meanwhile, explicit constructions remain quite elusive. The main way to triangulate $\RP^d$ is to construct a triangulation of $\bbS^d$ that forms a simplicial double cover of $\RP^d$. A classical result of K\"uhnel~\cite{vonKuhnel1987} shows that the barycentric subdivision of the boundary of the $(d+1)$-simplex forms a simplicial double cover of $\RP^d$, thereby providing a triangulation of $\RP^d$ with $2^{d+1}-1$ vertices. Recently, Venturello and Zheng~\cite{Venturello2021} improved this to $\Theta(\varphi^d)$ vertices where $\varphi\approx 1.618$ is the golden ratio. Soon after, Adiprasito, Avvakumov and Karasev~\cite{Adiprasito2022}  obtained a construction using at most $(\sqrt{d}+1)^{\sqrt{d}+1}\cdot 2^{\sqrt{d}}$ vertices. This breakthrough provides the first subexponential triangulations of the real projective space.
\\

\begin{figure}[h!]
  \centering
  \begin{minipage}[c]{0.25\textwidth}
    \centering
    \includegraphics[width=\linewidth]{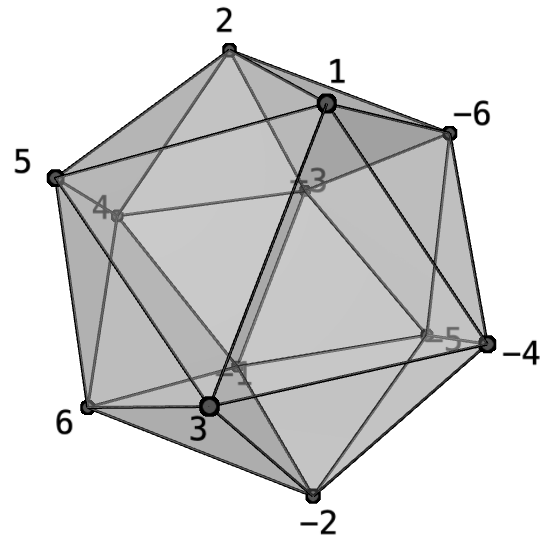}
  \end{minipage}
  \hspace{0.03\textwidth}
  \begin{minipage}[c]{0.06\textwidth}
    \centering
    {\Huge $\to$}
  \end{minipage}
  \hspace{0.03\textwidth}
  \begin{minipage}[c]{0.25\textwidth}
    \centering
    \includegraphics[width=\linewidth]{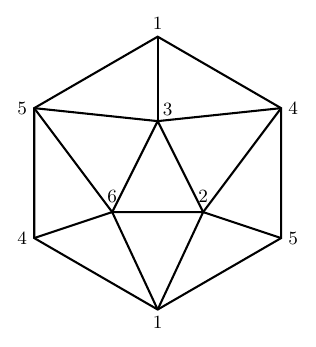}
  \end{minipage}

  \caption{Identifying antipodal vertices $\pm i$ for $i \in [6]$ of the icosahedron yields the unique minimal triangulation of $\RP^2$.}
  \label{fig:ico-to-rp2}
\end{figure}

Provably minimal triangulations of $\RP^d$ are only known in dimensions $1,2,3$, and $4$. For $d=1$, this is simply the (boundary of a) triangle as $\RP^1$ is homeomorphic to $\bbS^1$. For $d=2$, it is well-known that antipodally identifying the vertices of the icosahedron forms the unique minimal triangulation of $\RP^2$ (see Figure~\ref{fig:ico-to-rp2}). For $\RP^3$, Walkup~\cite{Walkup_1970} constructed a triangulation on $11$ vertices and proved that it is minimal. In low dimensions, computational methods have been effective; we specifically mention the program BISTELLAR \cite{Bjorner2000, Lutz2006} which uses bistellar flips, equivalently Pachner moves~\cite{Pachner1987,Pachner1991}. Since these methods are purely computational, the output is usually a list of facets without any geometric interpretation; this is what happens with Lutz's construction \cite{Lutz1999} of a minimal $16$-vertex triangulation of $\RP^4$. Later, its appropriate geometric interpretation was given by Balagopalan~\cite{Balagopalan}. As for $\RP^5$, the program also provides a list of $676$ facets that form a triangulation of $\RP^5$ with $24$ vertices~\cite{Lutz2006,BenedettiLutz:RP5_24} with little geometric explanation and no guarantee that this triangulation can be obtained via an antipodal quotient of the boundary of a polytope.

\subsection{Main result}
Our main result is a new $24$-vertex triangulation of $\RP^5$ which has a simplicial double cover that is the boundary of a convex $6$-polytope. This $6$-polytope is completely explicit, has simple coordinates with a combinatorially rich structure; it has a surprisingly  large number of symmetries as the automorphism group has order 192. We believe this polytope to be of independent interest.

 \begin{thm}[Main Result]
    There is a $48$-vertex centrally symmetric simplicial $6$-polytope for which the antipodal quotient of its boundary is a triangulation of $\RP^5$.  
\end{thm}

 We conjecture that this triangulation is vertex-minimal. We hope that this new geometric example will offer combinatorial insights into the triangulation problem for $\RP^d$ more broadly, introduce a new computational approach, and provide a highly symmetric $6$-polytope that may be of independent interest. The discovery of the polytope was computer-assisted and used black-box optimization in combination with an $L^1$-sparsity trick (see \S 3 for some more details). As for $\RP^6$, a construction of Venturello--Zheng \cite{Venturello2021} achieves a triangulation using 53 vertices. We are able to give a clean conceptual example using 49 vertices (using our main result) and, using the same computational approach, an explicit example using 45 vertices. This is potentially far from minimal as the Arnoux--Marin result gives a lower bound of $29$ vertices.

\section{Triangulations of $\RP^5$ and $\RP^6$}\label{sec:examples}
\subsection{Background} By a \emph{polytope} we mean the convex hull of finitely many points in a Euclidean space. Given a manifold $M$, a \emph{triangulation} of $M$ is a simplicial complex whose geometric realization is homeomorphic to $M$. A polytope $P$ is \emph{centrally symmetric} if $P = -P$. A simplicial complex is centrally symmetric if it has a free simplicial involution. For a simplicial complex $\Delta$ and a face $\sigma \in \Delta$, the (closed) \emph{star} of $\sigma$ in $\Delta$, denoted $\str_\Delta \sigma$, is the subcomplex of $\Delta$ generated by the facets of $\Delta$ containing $\sigma$.
Let $\Delta$ be a centrally symmetric simplicial $d$-sphere with a free involution $\tau$. Note that no vertex $v \in \Delta$ can be adjacent to its \emph{antipode} $\tau(v)$, otherwise $\tau$ would fix the edge $\{v, \tau(v)\}$. 
If \begin{equation}\label{eqn:disjointstar}
    \str_\Delta v \cap \str_\Delta\tau(v) = \{\emptyset\} \text{ for every vertex $v$ of $\Delta$},
\end{equation} then the complex $\Delta/\tau$ obtained by identifying every vertex $v$ with $\tau(v)$ is a triangulation of the $d$-dimensional real projective space $\RP^d$. This can be easily verified by using the fact that the $d$-sphere $\bbS^d$ is the double cover of $\RP^d$ and the definition of a simplicial complex. Condition~\eqref{eqn:disjointstar} is equivalent to the condition that no antipodal pair $v$ and $\tau(v)$ share a common neighbor in the $1$-skeleton of $\Delta$.

For further exposition on polytopes and triangulated manifolds, see \cite{Ziegler1995} and \cite{RourkeSanderson1982}, respectively.

\subsection{A $24$-vertex triangulation of $\RP^5$}
\begin{theorem} \label{thm:main}
For $\alpha=3/7, \beta=4/7$ and $\gamma=5/7$, the convex hull of the points

{\footnotesize \begin{longtable}{lll}
$P_{1} = (0, \gamma, \alpha, \beta, 0, 0)$ & $P_{2} = (\alpha, 0, 0, -\gamma, -\beta, 0)$ & $P_{3} = (\beta, \gamma, -\alpha, 0, 0, 0)$ \\
$P_{4} = (-\beta, \gamma, -\alpha, 0, 0, 0)$ & $P_{5} = (-\gamma, 0, 0, -\alpha, -\beta, 0)$ & $P_{6} = (0, 0, \beta, 0, \alpha, -\gamma)$ \\
$P_{7} = (0, 0, \beta, 0, -\gamma, -\alpha)$ & $P_{8} = (-\gamma, 0, 0, -\alpha, \beta, 0)$ & $P_{9} = (0, \beta, 0, 0, -\alpha, -\gamma)$ \\
$P_{10} = (0, -\gamma, -\alpha, \beta, 0, 0)$ & $P_{11} = (0, -\beta, 0, 0, \gamma, -\alpha)$ & $P_{12} = (0, -\beta, 0, 0, -\gamma, \alpha)$ \\
$P_{13} = (-\gamma, 0, 0, \alpha, 0, \beta)$ & $P_{14} = (-\gamma, 0, 0, \alpha, 0, -\beta)$ & $P_{15} = (0, 0, \beta, 0, -\alpha, \gamma)$ \\
$P_{16} = (\alpha, 0, 0, \gamma, 0, \beta)$ & $P_{17} = (0, \alpha, -\gamma, \beta, 0, 0)$ & $P_{18} = (\alpha, 0, 0, -\gamma, \beta, 0)$ \\
$P_{19} = (\beta, -\alpha, -\gamma, 0, 0, 0)$ & $P_{20} = (0, \beta, 0, 0, \alpha, \gamma)$ & $P_{21} = (0, 0, \beta, 0, \gamma, \alpha)$ \\
$P_{22} = (\beta, \alpha, \gamma, 0, 0, 0)$ & $P_{23} = (\alpha, 0, 0, \gamma, 0, -\beta)$ & $P_{24} = (0, \alpha, -\gamma, -\beta, 0, 0)$ \\
$P_{25} = (0, -\gamma, -\alpha, -\beta, 0, 0)$ & $P_{26} = (-\alpha, 0, 0, \gamma, \beta, 0)$ & $P_{27} = (-\beta, -\gamma, \alpha, 0, 0, 0)$ \\
$P_{28} = (\beta, -\gamma, \alpha, 0, 0, 0)$ & $P_{29} = (\gamma, 0, 0, \alpha, \beta, 0)$ & $P_{30} = (0, 0, -\beta, 0, -\alpha, \gamma)$ \\
$P_{31} = (0, 0, -\beta, 0, \gamma, \alpha)$ & $P_{32} = (\gamma, 0, 0, \alpha, -\beta, 0)$ & $P_{33} = (0, -\beta, 0, 0, \alpha, \gamma)$ \\
$P_{34} = (0, \gamma, \alpha, -\beta, 0, 0)$ & $P_{35} = (0, \beta, 0, 0, -\gamma, \alpha)$ & $P_{36} = (0, \beta, 0, 0, \gamma, -\alpha)$ \\
$P_{37} = (\gamma, 0, 0, -\alpha, 0, -\beta)$ & $P_{38} = (\gamma, 0, 0, -\alpha, 0, \beta)$ & $P_{39} = (0, 0, -\beta, 0, \alpha, -\gamma)$ \\
$P_{40} = (-\alpha, 0, 0, -\gamma, 0, -\beta)$ & $P_{41} = (0, -\alpha, \gamma, -\beta, 0, 0)$ & $P_{42} = (-\alpha, 0, 0, \gamma, -\beta, 0)$ \\
$P_{43} = (-\beta, \alpha, \gamma, 0, 0, 0)$ & $P_{44} = (0, -\beta, 0, 0, -\alpha, -\gamma)$ & $P_{45} = (0, 0, -\beta, 0, -\gamma, -\alpha)$ \\
$P_{46} = (-\beta, -\alpha, -\gamma, 0, 0, 0)$ & $P_{47} = (-\alpha, 0, 0, -\gamma, 0, \beta)$ & $P_{48} = (0, -\alpha, \gamma, \beta, 0, 0)$\\
\end{longtable}}
\noindent
is a centrally symmetric simplicial $6$-polytope such that the antipodal quotient of its boundary is a $24$-vertex triangulation of $\RP^5$. 
\end{theorem}

\noindent
\subsection*{A Simplicial Centrally Symmetric Polytope}  Let $\mathcal{P}_{6,48}$ denote this polytope. Because all of these vectors lie on the Euclidean $5$-sphere (of radius $\sqrt{50/49}$), we can conclude that each one of these points is a vertex of the polytope.  Additionally, the points are listed in a way that $P_i$ and $P_{i+24}$ are antipodal for all $i\leq 24$ which shows that the polytope is centrally symmetric. In particular, the boundary of $\mathcal{P}_{6,48}$ is homeomorphic to $\mathbb{S}^5$ and admits the free involution $x \mapsto -x$. Finally, with SageMath~\cite{sagemath}, we verify that $\mathcal{P}_{6,48}$ is simplicial and satisfies Condition~\eqref{eqn:disjointstar}. Furthermore, we discuss why Condition~\eqref{eqn:disjointstar} is satisfied when analyzing the $1$-skeleton.

\subsection*{Facets and Symmetries}
This polytope has $f$-vector 
$$(48, 552, 2432, 4776, 4272, 1424).$$ 
Division by $2$ gives the $f$-vector of the resulting triangulation of $\RP^5$: $$(24, 276, 1216, 2388, 2136, 712).$$

We note that the $f$-vector is different from that of the $24$-vertex triangulation obtained by Lutz, and therefore the triangulations are nonisomorphic. 
Roughly speaking, we can view $\mathcal{P}_{6,48}$ as a highly symmetrical system of overlapping cubes.  Specifically, let
\[
\mathcal{S} = \{\{1,2,3\},\{1,4,5\},\{1,4,6\},\{2,3,4\},\{2,5,6\},\{3,5,6\}\}.
\]
Then $V(\mathcal{P}_{6,48})$ can be partitioned into $V_S$ for $S \in \mathcal{S}$, where $V_S$ denotes the set of vertices whose support is exactly the coordinates indexed by $S$. See Table~\ref{tab:groups_of_cubes} below.

\begin{table}[H]
\footnotesize
\centering
\renewcommand{\arraystretch}{1.2}
\begin{tabular}{|c|c|c|}
\hline 
\begin{tabular}{l}
$\mathbf{S=\{1,4,5\}}$ \\
$\pm(\gamma,0,0,\alpha,\beta,0)$ \\
$\pm(\alpha,0,0,-\gamma,\beta,0)$ \\\
$\pm(-\gamma,0,0,-\alpha,\beta,0)$ \\
$\pm(-\alpha,0,0,\gamma,\beta,0)$ \\
\end{tabular}
&
\begin{tabular}{l}
$\mathbf{S=\{1,2,3\}}$ \\
$\pm(\beta,\alpha,\gamma,0,0,0)$ \\
$\pm(\beta,-\gamma,\alpha,0,0,0)$ \\
$\pm(\beta,-\alpha,-\gamma,0,0,0)$ \\
$\pm(\beta,\gamma,-\alpha,0,0,0)$ \\
\end{tabular}
&
\begin{tabular}{l}
$\mathbf{S=\{2,5,6\}}$ \\
$\pm(0,\beta,0,0,\alpha,\gamma)$ \\
$\pm(0,\beta,0,0,-\gamma,\alpha)$ \\
$\pm(0,\beta,0,0,-\alpha,-\gamma)$ \\
$\pm(0,\beta,0,0,\gamma,-\alpha)$ \\
\end{tabular}
\\ \hline 

\begin{tabular}{l}
$\mathbf{S=\{1,4,6\}}$ \\
$\pm(\alpha,0,0,\gamma,0,\beta)$ \\
$\pm(-\gamma,0,0,\alpha,0,\beta)$ \\
$\pm(-\alpha,0,0,-\gamma,0,\beta)$ \\
$\pm(\gamma,0,0,-\alpha,0,\beta)$ \\
\end{tabular}
&
\begin{tabular}{l}
$\mathbf{S=\{2,3,4\}}$ \\
$\pm(0,\gamma,\alpha,\beta,0,0)$ \\
$\pm(0,\alpha,-\gamma,\beta,0,0)$ \\
$\pm(0,-\gamma,-\alpha,\beta,0,0)$ \\
$\pm(0,-\alpha,\gamma,\beta,0,0)$ \\
\end{tabular}
&
\begin{tabular}{l}
$\mathbf{S=\{3,5,6\}}$ \\
$\pm(0,0,\beta,0,\gamma,\alpha)$ \\
$\pm(0,0,\beta,0,\alpha,-\gamma)$ \\
$\pm(0,0,\beta,0,-\gamma,-\alpha)$ \\
$\pm(0,0,\beta,0,-\alpha,\gamma)$ \\
\end{tabular}
\\ \hline 
\end{tabular}
\caption{Vertices of $\mathcal{P}_{6,48}$ grouped by their supports.}
\label{tab:groups_of_cubes}
\end{table}

Observe that $\conv(V_S)$ is a combinatorial $3$-cube for each $S$. The edges in each cube are given by pairs of vertices whose signs differ in exactly one coordinate. By examining the SageMath output of the list of facets of $\mathcal{P}_{6,48}$, we see that every facet is the convex hull (combinatorially, the \emph{join}) of faces from different cubes. Since the facets are simplicial, the maximal contribution from one cube is an edge.

Next, we show that $\mathcal{P}_{6,48}$ has an automorphism group of order $192$.
Let $Z_1, Z_2, Z_3$ respectively denote the three orthogonal planes spanned by the standard basis vectors $\{\e_1,\e_4\}, \{\e_2,\e_3\},\{\e_5,\e_6\}$ (these pairs correspond exactly to the shared nonzero coordinates of each column of Table~\ref{tab:groups_of_cubes}). Note that every linear automorphism of our point set will be an automorphism of $\mathcal{P}_{6,48}$. From Table~\ref{tab:groups_of_cubes}, we can see a dihedral symmetry within each of $Z_1, Z_2, Z_3$, a $3$-cycle permutation $Z_1Z_3Z_2$ of the planes, and the central symmetry. Consider two automorphisms
$$
b = \begin{pmatrix}
-1 &  0&  0&  0&  0& 0\\
 0& 1 &  0&  0&  0& 0\\
 0&  0& 1 &  0&  0& 0\\
 0&  0&  0& 1 &  0& 0\\
 0&  0&  0&  0& 0 & -1\\
 0&  0&  0&  0& 1 & 0
\end{pmatrix},
\qquad
c = \begin{pmatrix}
0 & 0 & 1 & 0 & 0 & 0\\
0 & 0 & 0 & 0 & 1 & 0\\
0 & 0 & 0 & 0 & 0 & -1\\
0 & -1 & 0 & 0 & 0 & 0\\
0 & 0 & 0 & -1 & 0 & 0\\
-1 & 0 & 0 & 0 & 0 & 0
\end{pmatrix}.
$$

Each matrix takes $\e_i$ to the $i$-th column of that matrix. Therefore, $b$ negates the first coordinate of $Z_1$, fixes $Z_2$, and applies a $90^\circ$ rotation to $Z_3$. Meanwhile, $c$ takes $\e_1, \e_4$ to $-\e_6, -\e_5$ (thus taking $Z_1$ to $Z_3$), $\e_2, \e_3$ to $-\e_4, \e_1$ (thus taking $Z_2$ to $Z_1$), and $\e_5, \e_6$ to $\e_2, -\e_3$ (thus taking $Z_3$ to $Z_2$). In particular, $c$ cyclically permutes the three planes with the $3$-cycle $Z_1Z_3Z_2$.

Let $G$ be the group generated by $b$ and $c$. We view $G$ as a group acting on the set of ordered triples 
\[
\{(Z_1, Z_2, Z_3), (Z_3, Z_1, Z_2), (Z_2, Z_3, Z_1)\}.
\]
Thus $|\operatorname{Orb}_G((Z_1, Z_2, Z_3))| = |\langle c \rangle| = 3$. Meanwhile, $\operatorname{Stab}_G((Z_1, Z_2, Z_3))$ must be generated by $b, cbc^{-1}, c^2bc^{-2}$, because whenever we permute the planes we need to undo the action. 

Next, we show that $\operatorname{Stab}_G((Z_1, Z_2, Z_3)) = \langle b, cbc^{-1}, c^2bc^{-2} \rangle$ has order $64$. 
View the actions of $b,cbc^{-1},c^2bc^{-2}$ on $Z_1,Z_2,Z_3$ as elements of the group $D_4 \times D_4 \times D_4$. Using $r$ to denote the rotation action and $s$ to denote the reflection action, we find that these three generators of $\operatorname{Stab}_G((Z_1, Z_2, Z_3))$ may be interpreted as $(s,e,r),(e,r,s)$ and $(r,s,e)$ respectively. Observe $(r,s,e)(r,s,e)=(r^2,e,e)$ and that the subgroup $H\coloneq \langle (r^2,e,e),(e,r^2,e),(e,e,r^2)\rangle$ is the commutator subgroup of $D_4 \times D_4 \times D_4$. The quotient group $\operatorname{Stab}_G((Z_1, Z_2, Z_3))/H$ can be expressed as
$$\frac{\langle (s,e,r),(r,s,e),(e,r,s)\rangle}{\langle(r^2,e,e),(e,r^2,e),(e,e,r^2)\rangle}.$$
Moreover, this quotient group is commutative and consists of eight distinct elements. Notice also that $|H|=2^3$. Hence, by Lagrange's theorem, we find that 
\begin{align*}
|\operatorname{Stab}_G((Z_1, Z_2, Z_3))|=
|\operatorname{Stab}_G((Z_1, Z_2, Z_3))/H||H|=2^3\cdot 2^3.
\end{align*}
Therefore, we can conclude that 
\begin{align*}
    |G| = |\operatorname{Orb}_G((Z_1, Z_2, Z_3))| \cdot |\operatorname{Stab}_G((Z_1, Z_2, Z_3))| = 3\cdot 2^6 = 192.
\end{align*}

Let $\partial \mathcal{P}_{6,48}$ denote the simplicial $5$-sphere that is the boundary of $\mathcal{P}_{6,48}$. By a SageMath computation, the order of $\operatorname{Aut}(\partial\mathcal{P}_{6,48})$ (all combinatorial automorphisms of the simplicial sphere), which contains $\operatorname{Aut}(\mathcal{P}_{6,48})$ (all geometric symmetries of the polytope), is $192$. Thus, we conclude that $G = \operatorname{Aut}(\mathcal{P}_{6,48}) = \operatorname{Aut}(\partial \mathcal{P}_{6,48})$.

\subsection*{1-skeleton} 
\begin{figure}[h!]
\centering
\includegraphics[width=10cm, height=8cm, keepaspectratio]{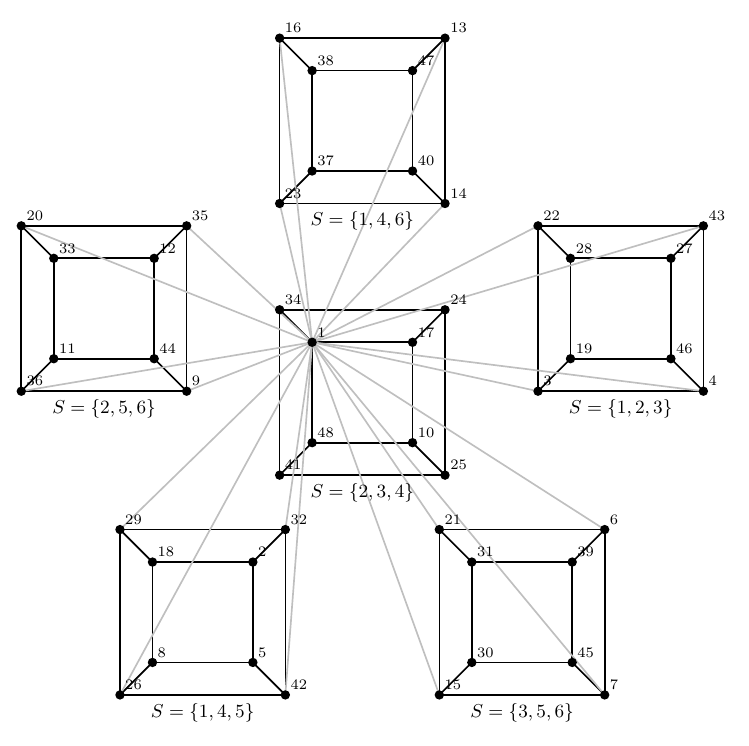}
\caption{The union of induced subgraphs on $V_S$ for every $S \in \mathcal{S}$ and the neighborhood of $P_1$.}
\label{fig:graph_nbhd}
\end{figure}

Observe that no two vertices of $\mathcal{P}_{6,48}$ have disjoint supports. Moreover, because of our choice of $\alpha, \beta, \gamma$, we can ensure that no two of these vectors have an inner product of zero. With SageMath, we can verify that there is an edge between vertices $p$ and $q$ if and only if $\langle p,q\rangle>0$. Consequently, the graph is $23$-regular on $48$ vertices and it satisfies Condition~\eqref{eqn:disjointstar}--- the condition that every vertex and its antipode have disjoint sets of neighbors. Therefore, because the graph is $23$-regular on $48$ vertices, it is tight with respect to this property. Furthermore, because the graph of $\mathcal{P}_{6,48}$ is $23$-regular, we remark that the graph of the triangulation of $\RP^5$ is isomorphic to $K_{24}$. 
Another description of a vertex neighborhood is as follows. Let $v \in V(\mathcal{P}_{6,48})$ with support $S(v) \in \mathcal{S}$. Then besides the three neighbors of $v$ in the cube $\conv(V_{S(v)})$, $v$ is connected to exactly one square from each of the other cubes. See Figure~\ref{fig:graph_nbhd} for an illustration.

We note that one may consider the graph $G_t=(V,E_t)$ with $V$ corresponding to the 48 vertices of $\mathcal{P}_{6,48}$ and the edge set
$ E_t = \left\{\{i,j\}: i \neq j ~ \mbox{and} ~ \left\langle p_i, p_j \right\rangle > t\right\}$ for some parameter $t \in \mathbb{R}$. The choice $t=19/49$ leads to a $10$-regular graph with $240$ edges and chromatic number $4$, $t= 17/49$ leads to an $11$-regular graph with $264$ edges and chromatic number $6$, $t= 15/49$ leads to a $15$-regular graph with $360$ edges with chromatic number $7$ and, finally, $t=11/49$ leads to a $23$-regular graph with $552$ edges and chromatic number $12$. We are not aware of any prior appearance of any of these graphs in the literature.

\subsection*{Links}
\begin{figure}[h!]
\centering
\includegraphics[width=7.5cm, height=5.5cm, keepaspectratio]{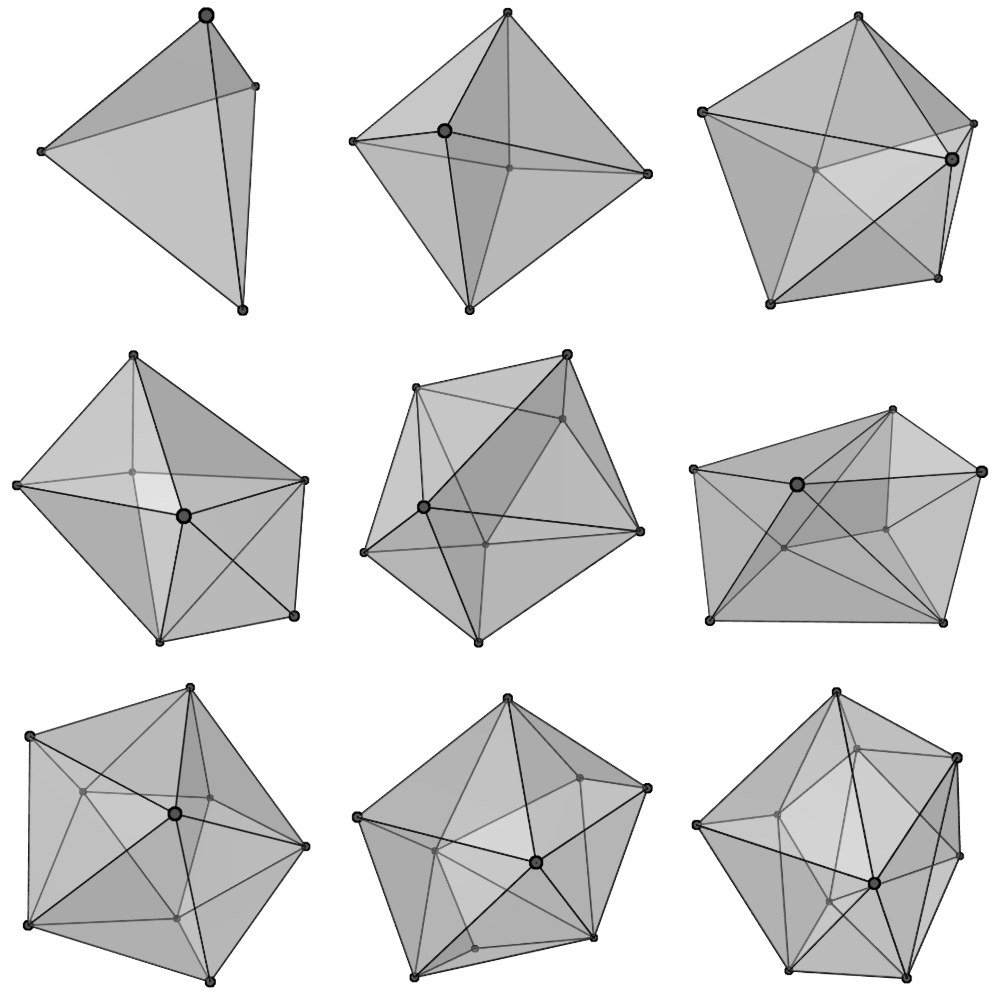}
\caption{Polytopal realizations of different combinatorial types of links of the $2$-faces of $\partial \mathcal{P}_{6,48}$, ordered by the number of vertices.}
\label{fig:triang_links}
\end{figure}

From a computation in SageMath, one finds that all vertex links of $\partial \mathcal{P}_{6,48}$ have the same combinatorial type, a certain $4$-sphere on $23$ vertices and $178$ facets. For $1 \le i \le 3$, the set of links of all $i$-faces contains more than one combinatorial type. In particular, the set of links of the $2$-faces contains $9$ different combinatorial types. Figure \ref{fig:triang_links} shows polytopal realizations of these $2$-spheres.

\subsection*{Comparison to the Lutz 
Triangulation} Recall that Lutz~\cite{Lutz1999} has already acquired a $24$-vertex triangulation of $\RP^5$; we take a moment to compare our triangulation with this one~\cite{BenedettiLutz:RP5_24}. It has an $f$-vector of $$(24,273,1174,2277,2028,676).$$ From this, one finds that it is three edges away from having a complete $1$-skeleton as $273=\binom{24}{2}-3$. Additionally, a computation in SageMath shows that the automorphism group of Lutz's triangulation has order $12$. In this sense, our construction is more symmetric.

\subsection{Triangulations of $\mathbb{R}P^6$}
We can easily obtain a $98$-vertex centrally symmetric simplicial $7$-polytope satisfying Condition (\ref{eqn:disjointstar}) from $\mathcal{P}_{6,48}$: first consider the cylinder given by $\mathcal{P}_{6,48} \times [-1, 1] \subseteq \mathbb{R}^7$, then introduce a new pair of antipodal vertices $\pm v \in \mathbb{R}^7$ above and below the cylinder to cone respectively over $\mathcal{P}_{6,48} \times \{1\}$ and $\mathcal{P}_{6,48} \times \{-1\}$ so that we glue two pyramids over $\mathcal{P}_{6,48}$ above and below the cylinder to obtain $\mathcal{P}_{7,98}$. We claim that all antipodal vertices of $\mathcal{P}_{7,98}$ are contained in disjoint subsets of facets. Indeed, this is true for $\mathcal{P}_{6,48}$, and this property is preserved for $\mathcal{P}_{6,48} \times [-1,1]$ as all facets are precisely of the form: $\mathcal{P}_{6,48} \times \{-1\}$, $\mathcal{P}_{6,48} \times \{1\}$ and $F \times [-1,1]$ where $F$ is a facet of $\mathcal{P}_{6,48}$. Furthermore, coning with $\pm v$ preserves this property. We may symmetrically perturb all pairs of vertices of $\mathcal{P}_{7,98}$. In doing so, we obtain a simplicial polytope whose boundary is a simplicial subdivision of the boundary of $\mathcal{P}_{7,98}$. The result forms a simplicial double cover of $\mathbb{R}P^6$.

With our computational methods, we can obtain a simplicial $7$-polytope, denoted $\mathcal{P}_{7,90}$, with $90$ vertices such that the antipodal quotient of its boundary is a $45$-vertex triangulation of $\RP^6$. This is the best currently known triangulation of $\RP^6$ in terms of minimizing the number of vertices, beating the prior record of $53$ vertices from Venturello--Zheng~\cite{Venturello2021} (see Remark~\ref{rem:Kolosov} for a further update). The vertices of $\mathcal{P}_{7,90}$ are listed below.

{\footnotesize
\begin{longtable}{ll}

$\pm \left(\frac{103}{134},0,0,0,0,-\frac{27}{73},\frac{107}{205}\right)$ &
$\pm \left(0,-\frac{25}{102},\frac{79}{99},\frac{16}{123},0,\frac{15}{104},-\frac{17}{33}\right)$ \\[0.5em]

$\pm \left(0,\frac{87}{175},\frac{2}{13},-\frac{74}{141},\frac{17}{50},\frac{31}{116},\frac{63}{122}\right)$ &
$\pm \left(\frac{46}{141},-\frac{73}{147},0,0,-\frac{47}{80},\frac{5}{74},\frac{103}{189}\right)$ \\[0.5em]

$\pm \left(-\frac{76}{111},-\frac{47}{200},0,\frac{173}{389},-\frac{16}{99},0,\frac{109}{217}\right)$ &
$\pm \left(\frac{69}{173},0,\frac{3}{49},0,\frac{36}{89},-\frac{41}{90},-\frac{43}{63}\right)$ \\[0.5em]

$\pm \left(-\frac{16}{181},\frac{19}{70},\frac{93}{137},-\frac{13}{61},-\frac{31}{61},-\frac{3}{49},\frac{19}{49}\right)$ &
$\pm \left(0,-\frac{59}{109},\frac{55}{274},-\frac{32}{91},\frac{67}{95},-\frac{28}{131},0\right)$ \\[0.5em]

$\pm \left(0,-\frac{6}{53},-\frac{16}{67},0,\frac{28}{71},\frac{38}{75},-\frac{77}{107}\right)$ &
$\pm \left(\frac{31}{69},\frac{81}{121},0,0,-\frac{7}{85},0,\frac{58}{99}\right)$ \\[0.5em]

$\pm \left(0,\frac{24}{83},\frac{24}{43},-\frac{129}{212},\frac{31}{92},-\frac{8}{23},0\right)$ &
$\pm \left(-\frac{9}{35},0,-\frac{59}{100},\frac{49}{100},0,-\frac{1}{4},-\frac{33}{62}\right)$ \\[0.5em]

$\pm \left(-\frac{143}{199},0,\frac{55}{102},-\frac{47}{107},0,0,0\right)$ &
$\pm \left(\frac{7}{124},-\frac{56}{337},-\frac{21}{79},0,-\frac{44}{75},\frac{66}{115},-\frac{46}{97}\right)$ \\[0.5em]

$\pm \left(\frac{3}{26},0,\frac{51}{148},\frac{45}{92},\frac{19}{39},-\frac{18}{65},-\frac{23}{41}\right)$ &
$\pm \left(\frac{15}{73},\frac{49}{60},-\frac{9}{29},-\frac{18}{139},\frac{13}{34},-\frac{7}{52},-\frac{5}{43}\right)$ \\[0.5em]

$\pm \left(\frac{103}{205},0,\frac{10}{37},0,-\frac{107}{140},\frac{43}{143},0\right)$ &
$\pm \left(\frac{10}{97},-\frac{34}{77},-\frac{6}{41},0,\frac{29}{79},\frac{21}{83},\frac{97}{128}\right)$ \\[0.5em]

$\pm \left(-\frac{83}{145},\frac{32}{45},\frac{10}{33},\frac{14}{61},-\frac{5}{49},0,-\frac{6}{55}\right)$ &
$\pm \left(-\frac{17}{61},-\frac{8}{67},-\frac{32}{55},0,\frac{46}{89},\frac{11}{20},0\right)$ \\[0.5em]

$\pm \left(-\frac{9}{34},-\frac{22}{35},-\frac{29}{40},-\frac{11}{114},0,0,0\right)$ &
$\pm \left(\frac{21}{58},0,-\frac{25}{38},0,-\frac{55}{108},-\frac{23}{57},-\frac{17}{145}\right)$ \\[0.5em]

$\pm \left(\frac{19}{28},-\frac{41}{90},-\frac{5}{77},0,0,0,-\frac{67}{117}\right)$ &
$\pm \left(-\frac{27}{46},-\frac{63}{134},0,-\frac{40}{87},\frac{8}{93},0,\frac{13}{28}\right)$ \\[0.5em]

$\pm \left(0,\frac{12}{91},-\frac{65}{86},-\frac{49}{137},0,0,-\frac{41}{77}\right)$ &
$\pm \left(-\frac{31}{75},0,-\frac{61}{203},0,-\frac{83}{110},-\frac{7}{17},0\right)$ \\[0.5em]

$\pm \left(-\frac{58}{75},0,-\frac{35}{132},0,0,\frac{69}{121},\frac{5}{63}\right)$ &
$\pm \left(-\frac{40}{213},-\frac{66}{131},0,-\frac{39}{112},0,-\frac{18}{41},-\frac{46}{73}\right)$ \\[0.5em]

$\pm \left(0,-\frac{21}{100},0,-\frac{33}{74},\frac{107}{135},\frac{27}{88},-\frac{21}{113}\right)$ &
$\pm \left(0,-\frac{8}{39},0,-\frac{18}{61},-\frac{9}{112},-\frac{97}{136},\frac{65}{109}\right)$ \\[0.5em]

$\pm \left(\frac{17}{106},0,-\frac{31}{103},-\frac{31}{202},0,-\frac{91}{103},-\frac{35}{124}\right)$ &
$\pm \left(\frac{68}{131},\frac{51}{142},-\frac{17}{84},-\frac{44}{81},-\frac{705}{1411},0,\frac{20}{159}\right)$ \\[0.5em]

$\pm \left(0,-\frac{47}{182},\frac{16}{127},-\frac{64}{91},-\frac{3}{26},\frac{16}{25},0\right)$ &
$\pm \left(0,\frac{97}{135},\frac{41}{101},0,-\frac{33}{127},\frac{175}{349},0\right)$ \\[0.5em]

$\pm \left(-\frac{79}{201},0,\frac{13}{43},\frac{77}{108},-\frac{41}{87},\frac{51}{331},0\right)$ &
$\pm \left(-\frac{14}{43},\frac{301}{502},-\frac{74}{117},\frac{21}{131},-\frac{6}{47},\frac{31}{102},0\right)$ \\[0.5em]

$\pm \left(-\frac{14}{99},-\frac{31}{39},\frac{3}{34},0,0,\frac{7}{12},0\right)$ &
$\pm \left(0,-\frac{59}{83},0,\frac{38}{99},0,-\frac{33}{56},0\right)$ \\[0.5em]

$\pm \left(\frac{11}{51},-\frac{41}{85},\frac{11}{84},-\frac{589}{1177},-\frac{52}{97},-\frac{64}{183},-\frac{19}{91}\right)$ &
$\pm \left(-\frac{13}{150},-\frac{15}{71},-\frac{35}{83},\frac{19}{31},0,0,\frac{76}{121}\right)$ \\[0.5em]

$\pm \left(-\frac{223}{322},\frac{23}{134},-\frac{12}{23},\frac{23}{242},-\frac{3}{32},-\frac{13}{29},0\right)$ &
$\pm \left(\frac{43}{104},-\frac{23}{102},\frac{357}{1072},\frac{2}{3},-\frac{11}{61},\frac{17}{48},\frac{29}{114}\right)$ \\[0.5em]

$\pm \left(-\frac{46}{91},0,\frac{31}{74},0,0,-\frac{80}{109},-\frac{15}{86}\right)$ &
$\pm \left(\frac{27}{92},-\frac{76}{157},0,\frac{129}{187},\frac{37}{83},\frac{24}{337},0\right)$ \\[0.5em]

$\pm \left(-\frac{13}{25},0,0,-\frac{37}{314},-\frac{46}{59},\frac{23}{120},-\frac{99}{371}\right)$ &
\\

\end{longtable}
}

The complete code verifying the properties of $\mathcal{P}_{6,48}$, $\mathcal{P}_{7,90}$, and Lutz's $24$-vertex triangulation of $\RP^5$ can be found in the notebook \nolinkurl{RP5_and_RP6_verification_code.ipynb} in \cite{ourcode}. The notebook runs under the SageMath 10.3 kernel.

\section{Remarks and Future Directions}\label{sec:future}

\subsection{Construction of $\mathcal{P}_{6,48}$} 
The construction of the example in $\mathbb{R}P^5$ was motivated by the following simple observation: suppose $P$ is a centrally symmetric simplicial polytope with the property that for any edge $\{x,y\}$ we have $\langle x,y \rangle >0$. Then there cannot be a vertex $w\in V(P)$ that is adjacent to a pair of antipodal vertices of $P$. This suggests an optimization problem: among all centrally symmetric sets of $48$ points on $\mathbb{S}^5$, find a configuration that maximizes the minimal inner product between any two points that are connected by an edge in the $1$-skeleton.  Geometrically, arrange the points in such a way that the maximal length of an edge of the convex hull is minimized (with the goal of it being $<\sqrt{2}$). This is an optimization problem in $5 \text{ coordinates} \times 48 \text{ points} = 240$ variables.
We used Google DeepMind's \textsc{AlphaEvolve}~\cite{novikov2025} as a way to do black-box optimization and ran a large number of instances in the range $44 \leq N \leq 64$ before ultimately finding a candidate set of 48 elements $x_1, \dots, x_{48} \in \mathbb{S}^5$ with minimal inner product $\sim 0.1$. Since the entire ansatz was invariant under orthogonal rotations, that set appeared to be entirely without structure. It was our hope that this set would be a rotation of a highly structured object containing a lot of 0 entries. This suggested looking at a real-valued map on the space of orthogonal matrices $f: Q_6 \rightarrow \mathbb{R}$ given by
$$ f(Q) = \sum_{i=1}^{48} \| Q x_i \|_{\ell^1} \qquad \mbox{and solving} \qquad f(Q) \rightarrow \min.$$
This exploits the classical idea that the $\ell^1$-norm promotes sparsity in its representation.  Approximately solving this problem led to a representation that allowed us to guess the structure to be the one described in the Main Theorem. We were then able to verify this with SageMath. It is noteworthy that 
the $2$-parameter family thus recovered contains examples with a lower bound on the inner product of $0.25$ which we were not able to find using the pure optimization approach.  

\subsection{Future Directions.}
Our main conjecture is that the triangulation of $\RP^5$ coming from $\mathcal{P}_{6,48}$ is vertex-minimal. We expect it to be a vertex-minimal triangulation of $\RP^5$, or, at the very least, to be vertex-minimal with respect to all triangulations coming from antipodally identifying a centrally symmetric polytope.
\begin{conjecture}
    The triangulation of $\RP^5$ coming from antipodally identifying $\mathcal{P}_{6,48}$ is vertex-minimal. 
\end{conjecture}
We have a few reasons to believe this. First, neither our approach, nor BISTELLAR~\cite{Lutz2006} could find a triangulation on fewer vertices. Second, the $1$-skeleton of our triangulation is tight with respect to the property described in Equation~\eqref{eqn:disjointstar}. That is, for every vertex $w$ and any pair of antipodal vertices $v,-v$, we have that $w$ is adjacent to exactly one of $v,-v$. Finally, we feel that any vertex-minimal triangulation of $\RP^d$ must come from a triangulated $d$-sphere that exhibits high degrees of symmetry. The fact that $48$ is highly composite may allow $\mathcal{P}_{6,48}$ to exhibit symmetries that any triangulation of $\mathbb{S}^5$ on $44$ or $46$ vertices cannot. 

As for the triangulation of $\RP^6$ coming from $\mathcal{P}_{7,90}$, the situation becomes more unclear. On one hand, our search could not find a $7$-polytope on fewer vertices whose quotient was $\RP^6$. On the other hand, this leaves the bounds for $\RP^6$ to lie between $29$ and $45$, which is quite large. Moreover, the K\"uhnel triangulations on $2^{d+1}-1$
vertices were shown to be quite suboptimal by the results of ~\cite{Basak_2013,Venturello2021,Adiprasito2022}. Therefore, one should not expect an increase by a factor of $2$ when moving from dimension $d$ to dimension $d+1$ for general $d$. 

\begin{question}\label{q:rp6minimality}
    Can one find a triangulation of $\RP^6$ on fewer than $45$ vertices?
\end{question}

\begin{remark}\label{rem:Kolosov}
    Since the publication of this article, Kolosov~\cite{Kolosov} has answered Question~\ref{q:rp6minimality} in the affirmative by producing a $44$-vertex triangulation of $\mathbb{R}P^6$. It remains intriguing to find the exact number of vertices needed to minimally triangulate $\mathbb{R}P^6$.
\end{remark}

\begin{remark}\label{rem:frick}
    Also since the publication of this article, the preprint of Frick et al.~\cite{frick2026} has significantly improved the lower bound of Arnoux--Marin~\cite{Arnoux_1991} by producing a superpolynomial lower bound on the number of vertices needed to triangulate $\mathbb{R}P^d$. 
\end{remark}

\textbf{Acknowledgements.}
The authors thank the anonymous referees for their comments to improve the exposition and rigor of this article. DG and YY thank Isabella Novik for her exemplary mentorship and fruitful discussions related to this problem, and Joshua Hinman for illustrative insights. Additionally, DG thanks Lorenzo Venturello for suggesting this problem. The authors are grateful to Google DeepMind for providing access to \textsc{AlphaEvolve} with special thanks to Bogdan Georgiev and Adam Zsolt Wagner.

\bibliographystyle{alpha}
\bibliography{Paperbiblio}
\end{document}